\input amstex%
\documentstyle{amsppt}%
\topmatter%

\title%
On perpetuities related to the size-biased distributions%
\endtitle%

\author%
Aleksander M. Iksanov%
\endauthor%

\address%
Cybernetics Faculty, Kiev National University,
01033 Kiev, Ukraine.
\endaddress%

\email%
iksan\@unicyb.kiev.ua.%
\endemail%


\abstract%
We study perpetuities of a special type related to the size-biased
distributions. Necessary and sufficient conditions of their existence and
uniqueness are obtained. A crucial point in proving all results is a close
connection between perpetuities treated in the paper and fixed points of so-called Poisson shot noise transforms.%
\endabstract%

\endtopmatter%

\document%

\head Introduction.\endhead%

Let $\Cal P^+$ be the set of all probability measures on the Borel subsets of $\Bbb R^+:=[0,\infty)$ and, for fixed $m>0$, $\Cal P_{m}^{+}:=\{\nu \in \Cal P^{+}:\int_{0}^{\infty }x\nu (dx)=m\}$. Recall that given distribution $\nu \in \Cal P_{m}^{+}$ , a distribution $\nu_{sb}$ is said to be the size-biased distribution corresponding to $\nu $, if
$$%
\nu _{sb}(dx)=m^{-1}x\nu (dx).%
$$%
Throughout this paper the symbols $\Cal L(\cdot )$ and $\Cal L%
_{sb}(\cdot )$ stand for the probability distribution and the size-biased
distribution of a random variable (rv) in question, respectively.

Consider the distributional equality
$$%
\ X \overset d \to = AX+B,%
\tag 1
$$%
where a random pair $(A,B)$ is independent of an rv $X$,
and ''$\overset d \to = $'' means ''equality of distributions''. In the recent
literature it is customary to say that the rv $X$ is a perpetuity. Cf. Embrechts, Goldie (1994) for more details.

In this paper we treat a quite special type of equality (1) in which $%
\Cal L(B)=\mu $, $\Cal L(X)=\mu _{sb}$, and $A$, $B$ are
independent rv's. It is convenient to put $B\equiv \eta $, $X\equiv \eta
_{sb}$ and rewrite (1) as follows
$$%
\eta _{sb} \overset d \to =A\eta _{sb}+\eta.%
\tag 2
$$%
Pitman, Yor (2000, p.35) mention the following problem: ''given a
distribution of $A$...whether there exists such a distribution of $\eta $''.

In all Propositions stated below we assume that $\Bbb P (A=0)=0$ and
keep the following notations $\mu :=\Cal L(\eta )$, $\rho :=\Cal L(A)$ where $\eta $ and $A$ satisfy (2). First we provide necessary and
sufficient conditions of the existence of \it {non-zero} \rm $\mu $ and reveal that when exists this distribution is unique up to the scale.

\proclaim {Proposition 1.1 (Existence and uniqueness)}%
\newline Assume that $\Bbb E \log A $ exists, finite or infinite. Then non-zero $\mu $ exists iff
$$%
\Bbb E \log A<0.%
\tag 3
$$%
Given $m>0$ there exists the unique distribution $\mu $ of mean $m$.%
\endproclaim%

In fact the above Proposition is essentially based upon a close relation
between solutions to (2) and fixed points of so-called Poisson
shot-noise transform (see Section 2 for definition and some properties).
This is the content of Proposition 2.1.

Given $\rho =\Cal L(A)\in \Cal P^{+}$ define a random map $\Bbb U_{\rho }:\Cal P_{m}^{+}\rightarrow \Cal P_{m}^{+}$ as follows:
$$%
\Cal L_{sb}(\Bbb U_{\rho }\theta )=\Bbb U_{\rho }\theta \ast
\Cal L(AX_{sb}),%
$$%
where $\theta =\Cal L(X)\in \Cal P_{m}^{+}$, $\Cal L%
(X_{sb})=\theta _{sb}$, $A$ and $X_{sb}$ are independent rv's and ''$\ast $%
''stands for the convolution of measures. To be sure that this map is
well-defined it suffices to verify that $\Bbb U_{\rho }\theta $ is of
finite mean. To this end consider the measure $N$ given by the equality $xN(dx):=\Cal L(AX_{sb})(dx)$. Clearly, it is the L\'{e}vy measure of some infinitely divisible (ID) distribution $\kappa $, say. Since $\int_{0}^{\infty }xN(dx)<\infty $ then $\int_{0}^{\infty }x\kappa(dx)<\infty $, and now it is obvious that $\kappa =\Bbb U_{\rho }\theta $.

In the proof of Proposition 1.1 it will be shown that $\mu $
is a weak limit of iterates $\Bbb U_{\rho}^{n}=\Bbb U_{\rho }(%
\Bbb U_{\rho }^{n-1})$, $n\in \Bbb N$, and therefore is a fixed point of $\Bbb U_{\rho }$. As it turned out, by putting to $\rho $ some additional moment
restrictions one may prove that $\Bbb U_{\rho }$ is a strict contraction
acting on some complete metric space ($T,r$). Using the Banach Fixed
Point Theorem allows us to assert that the iterates of $\Bbb U_{\rho }\theta $ converge exponentially fast to a unique fixed point of $\Bbb U_{\rho }$, for every $\theta \in T$.

For fixed $1<\Delta <2$, $m>0$ consider the set of probability
measures
$$%
{P}_{m}^{+}(\Delta ):=\{\nu \in P_{m}^{+}:\int_{0}^{\infty
}x^{\Delta }\nu (dx)<\infty \}.%
$$%
In view of Lemma 3.1 by Baringhaus, Gr\"{u}bel (1997) the quantity
$$%
\gather%
\ r_{\Delta }:=r_{\Delta }(\nu_{1},\nu_{2})=%
\\%
=\int_{0}^{\infty }s^{-\Delta -1}\left| \int_{0}^{\infty }\exp (isx)\nu
_{1}(dx)-\int_{0}^{\infty }\exp (isx)\nu _{2}(dx)\right| ds,%
\endgather%
$$%
defined for all $\nu_{1},\nu_{2}\in \Cal P_{m}^{+}(\Delta )$, is a metric
on $\Cal P_{m}^{+}(\Delta )$, and ($\Cal P_{m}^{+}(\Delta
),r_{\Delta })$ is a complete metric space.

In the next Proposition we study the case of $\rho $ having some finite
moments, so the function $g(x):=\Bbb EA^{x}$ is finite
(and log-convex) at least on some neighbourhood of the origin.

\proclaim {Proposition 1.2 (Contraction properties and results on moments)}%
a) Let
$$%
\text{ }g(p)<1\text{, for some }p>0\text{ },%
\tag 4
$$%
and $q\in (1,2)$ be any fixed number such that $g(q-1)<1$. Then \newline
1) $\Bbb U_{\rho }$ defined on a complete metric space $(\Cal P %
_{m}^{+}(q),r_{q})$ is a strict contraction, and therefore for every $\theta
\in \Cal P_{m}^{+}(q)$ the sequence $\Bbb U_{\rho }^{n}\theta $, $%
n=1,2,...$ converges in $r_{q}$-metric (hence, weakly) at an exponential
rate to a unique fixed point $\mu $ of $\Bbb U_{\rho }$;\newline
2) the rv $\eta $ with $\Cal L(\eta )=\mu $ solves (2), $\Bbb E\eta
=m$ and $\Bbb E\eta ^{1+p}<\infty $.\newline
b) If there exists non-zero $\eta $ satisfying (2) and $\Bbb E\eta ^{p+1}<\infty $%
, for some $p>0$, then $\Bbb E A^{p}<1$.%
\endproclaim%

It is easily seen that all $\mu $'s are ID. This follows from the
representation of non-negative ID distributions due to Steutel (see Sato
(1999, Theorem 51.1) for a detailed proof). In the next assertion we study
the structure of $\mu $ as an ID distribution more carefully.
\proclaim {Proposition 1.3 (Infinite divisibility)}%
All non-zero $\mu $'s are ID with the shift $0$ and the L\'{e}vy
measure $M(dx)=x^{-1}\Cal L(A\eta _{sb})(dx)$.
Furthermore, $\mu $'s are compound Poisson provided $x^{-1}\rho (dx)$ is
integrable at the neighbourhood of zero.%
\endproclaim%

Two other results of the paper deal with the tail behaviour of $\mu $. This
is changed considerably according to whether the rv $A$ may take values greater
than $1$, or not; if not, then whether esssup $\rho =1$ or not.
\proclaim {Proposition 1.4 (Exponential moments)}%
\newline a) $\mu $'s have finite exponential moment iff esssup $%
\rho \leq 1$; \newline
b) if esssup $\rho <1$ then $\mu $'s have entire characteristic functions.%
\endproclaim%
We also give a very simple independent proof of the next result one
implication of which follows by Proposition 1.4(a). Corollary 4.2 of Goldie, Gr\"{u}bel (1996) contains a result concerning general perpetuities in the spirit of the converse part of Proposition 1.5.

\proclaim {Proposition 1.5}%
Condition ''$\rho $ is concentrated on $(0$,$1]$'' is necessary and
sufficient to ensure the existence of the solution $\eta $ to (2) which is
completely determined by its moments.%
\endproclaim%
\bf{Notations and convention.} \rm''LT''(''LST'')-Laplace (Stieltjes)
transform, ''rv''-random variable, ''ID''- infinitely divisible,
''a.s.''-almost sure(ly), ''w.l.o.g.'' - without loss of generality; we always take distribution functions to be right-continuous.

\head Connection to fixed points of shot noise transforms.\endhead%

Throughout this Section we assume that all rv's involved live on a common
probability space $(\Omega $, $\Cal F$, $\Bbb P)$. Let $\xi $%
, $\xi _{1},\xi _{2},...$ be non-negative iidrv's, independent of the Poisson flow $\{\tau _{i}\},i \in \Bbb N$ with the intensity $0<\lambda <\infty $. Given $\Cal L(\xi _{i})$ fix a Borel measurable function $h:(0,\infty)\rightarrow \lbrack 0,\infty )$ enjoying the property $\int_{0}^{\infty }\Bbb E [1\wedge h(s)\xi ]ds<\infty $. Recall that under the above assumptions random series $\sum_{i=1}^{\infty }\xi _{i}h(\tau _{i})$
converges a.s., and its distribution is called a \it{(Poisson) shot noise
distribution.} \rm The function $h$ is said to be the \it response \rm function. We refer to Vervaat (1979) and Bondesson (1992, Section 10) for some additional
information regarding the shot noise distributions.

For a fixed $\lambda $, consider a \it {Poisson shot noise transform} \rm (SNT,
in short) $\Bbb T_{h,\lambda }$ acting on the set (domain)
$$%
\Cal P_{h}^{+}:=\{\nu \in \Cal P^{+}:\int_{0}^{\infty
}\int_{0}^{\infty }[1\wedge h(s)y]ds\nu (dy)<\infty \}%
$$%
with values in $\Cal P^{+}$ as follows
$$%
\Bbb T_{h,\lambda }(\Cal L(\xi )):=\Cal L\left(%
\sum_{i=1}^{\infty }\xi _{i}h(\tau _{i})\right).%
$$%

Thus the domain $\Cal P_{h}^{+}$ is the set of possible distributions
for rv $\xi $ that would ensure the well-definedness of shot noise
distribution. By \it {fixed points} \rm of the SNT we mean \it {non-zero} \rm
distributional solutions to the equation
$$%
\mu ^{\ast}=\Bbb T_{h,\lambda }(\mu ^{\ast}),%
$$%
where $\mu ^{\ast } =\Cal L(\xi )$. We will essentially make use of an
equivalent definition of fixed points given via the LST $\varphi
(s)=\int_{0}^{\infty }e^{-sx}\mu ^{\ast}(dx)$. Namely $\mu ^{\ast}\ $is a fixed point of
the SNT $\Bbb T_{h,\lambda }$ iff
$$%
\varphi (s)=\exp \left( \lambda \int_{0}^{\infty }(\varphi
(sh(u))-1)du\right).%
\tag 5
$$%
In what follows it is assumed that
$$%
\text{the response function }h\text{ is right-continuous and non-increasing.}
\tag 6
$$%
This assumption permits to define the right-continuous and non-increasing
\it {generalized inverse} \rm of $h$ given as follows $h^{\leftarrow }(z)=\inf
\{u:h(u)<z\}$, for $z<h(0^{+})$, and $0$, otherwise. Note that fixed points of
the SNT for $h$'s concentrated on $[0,1]$ and satisfying (6) have been
studied in Iksanov, Jurek (2002). They obtained necessary and sufficient
conditions of their existence and uniqueness.

The key ingredient for the proof of all assertions in Section 1 sounds as
follows

\proclaim {Proposition 2.1}%
a) If $\eta $ is an rv satisfying (2) with $\rho =\Cal L (A)$
concentrated on $[a,b]$ and such that $\Bbb E\log A<0$ then $\mu =%
\Cal L (\eta )$ is a fixed point of the SNT $\Bbb T_{h,1}$ \newline ($\mu
\equiv \mu ^{\ast }$) with $h$ defined via its generalized inverse
$$%
h^{\leftarrow }(x)=\int_{x}^{b}z^{-1}\rho (dz),x\in (a,b).%
\tag 7
$$%
Therefore $h$ is right-continuous and non-increasing with $\int_{0}^{\infty
}h(z)dz=1$ and \newline $\int_{0}^{\infty }h(z)\log h(z)dz<0$. \newline
b) Conversely, if $\mu ^{\ast }$ is a fixed point of the SNT $\Bbb T%
_{h,\lambda }$ with $h$ satisfying (6) and
$$%
\lambda \int_{0}^{\infty }h(z)dz=1\text{, }\int_{0}^{\infty }h(z)\log
h(z)dz<0\text{ and }\underset{z\rightarrow +0} \to \lim h(z)=b\in (0,\infty ]%
\,%
\tag 8
$$%
then an rv $\eta $ such that $\Cal L (\eta )=\mu ^{\ast }$ verifies (2)
with an rv $A$ whose distribution $\rho (dx)=-\lambda xh^{\leftarrow }(dx)$ is concentrated on $[a,b]$.%
\endproclaim%

\head The proofs.\endhead%

The first Lemma of this Section is implicit in Athreya (1969, Theorem 1).

\proclaim { Lemma 3.1}%
Let $\varphi _{1}(s)$ and $\varphi _{2}(s)$ be the LST's
of probability measures of the same (finite) mean. If for large enough
positive integers $n$ the function $\psi (s)=\dfrac{\left| \varphi
_{1}(s)-\varphi _{2}(s)\right| }{s}$ satisfies the inequality
$ \psi (s)\leq \Bbb E\psi (C_{n}s)\ ,$
where $\{C_{n}\}$ is a sequence of rv's tending to $0$ a.s.,
then $\varphi _{1}(s)\equiv $ $\varphi _{2}(s)$.%
\endproclaim%

The next Lemma is both an existence and uniqueness result concerning
fixed points of the SNT. In fact, the main part of Proposition 1.1 is merely a combination of Lemma 3.2 and Proposition 2.1.

\proclaim {Lemma 3.2}%
Let $h$ satisfies (6), $\lambda \int_{0}^{\infty
}h(z)dz=1$ and $\int_{0}^{\infty }h(z)\log h(z)dz<0$. Then given $m\in
(0,\infty )$ $\Bbb T_{h,\lambda }$ has a unique fixed point $\mu ^{\ast
} $ with $\int_{0}^{\infty }x\mu ^{\ast }(dx)=m$.%
\endproclaim%
\demo {Proof}%
For fixed $m>0$ consider the set of probability measures $%
\Cal P_{h,m}^{+}:=\{v\in \Cal P_{h}^{+}:\int_{0}^{\infty }xv(dx)=m\}
$. Starting with $\mu _{0}=\delta _{m}$, construct the sequence
$$%
\mu _{n}:=\Bbb T_{h,\lambda }\mu _{n-1}:=\Bbb T_{h,\lambda }^{n}\mu
_{0},n=1,2,...%
$$%
which is trivially well-defined on $\Cal P_{h,m}^{+}$ provided $%
\int_{0}^{\infty }h(z)dz<\infty $. The corresponding LST's $\varphi
_{n}^{(1)}(s)=\int_{0}^{\infty }e^{-sx}\mu _{n}(dx)$, $n=0,1,...$ satisfy
equations
$$%
\varphi _{0}^{(1)}(s)=e^{-ms}\ , \varphi _{n}^{(1)}(s)=\exp
\{-\lambda \int_{0}^{\infty }(1-\varphi _{n-1}(sh(u)))du\}\ , n=1,2,...%
\tag 9
$$%
Similarly, for fixed $s_{0}>0$, let us define yet another sequence $%
\{\varphi _{n}^{(2)}\}$ that satisfies (9) for $n=1,2,...$, but $1-\varphi
_{0}^{(2)}(s)=(1-e^{-s})1_{\{s\in \lbrack 0,s_{0}]\}}$, where $1_{A}$ is the
indicator of set $A$. \newline
Let us verify that the weak limit of $\mu _{n}$, as $n\rightarrow \infty $,
exists and has mean $m$. As it is well-known, this will mean that $\Bbb T%
_{h,\lambda }$ has a fixed point (which is this weak limit) on $\Cal P%
_{h,m}^{+}$.\newline
In what follows we use some ideas of Durrett, Liggett (1983, the proof of
Theorem 2.7). From (9) one gets, $\varphi _{1}^{(i)}(s)\geq $ $\varphi
_{0}^{(i)}(s)$, $i=1,2$  that implies
$$%
\varphi _{n}^{(i)}(s)\geq \varphi _{n-1}^{(i)}(s),\text{ }n=1,2,...,\text{ }%
s\geq 0\ , i=1,2.%
$$%
Thus the monotone and bounded sequence $\{\varphi _{n}^{(i)}\},n=1,2,...$
has a unique limit $\varphi ^{(i)}$, $i=1,2$, with
$$%
\underset{s\rightarrow +0} \to {\lim \sup } s^{-1}(1-\varphi ^{(i)}(s))\leq m\text{%
, }i=1,2.%
\tag 10
$$%
For what follows it is essential that $\varphi _{0}^{(2)}(s)\geq $ $\varphi
_{0}^{(1)}(s)$ implies
$$%
s^{-1}(1-\varphi ^{(2)}(s))\leq s^{-1}(1-\varphi ^{(1)}(s))\text{, }s\geq 0%
.%
\tag 11
$$%
Note $\varphi ^{(1)}(s)$ is the LST of a probability measure $\mu^{(1)}$, say.
Since $\int_{0}^{\infty }h(z)dz<\infty $ then by dominated convergence it is
easily seen that $\varphi ^{(1)}(s)$ satisfies the fixed point equation (5)
or equivalently $\mu^{(1)} $ is a (possibly degenerate at $0$) fixed point of the SNT. It remains to check that $\mu ^{(1)}\in \Cal P_{h,m}^{+}$.

To this end for $n=0,1,...$ put $\Phi _{n}(s):=e^{s}(1-\varphi
_{n}^{(2)}(e^{-s}))$, \newline $\Psi _{n}(s):=-e^{s}(1-\varphi
_{n}^{(2)}(e^{-s})+\log \varphi _{n}^{(2)}(e^{-s}))$. In view of assumptions
of the Lemma $\pi (dz):=-\lambda zh^{\leftarrow }(dz)$ is a probability
distribution. Let $\theta $, $\theta _{1}$, $\theta _{2}$, ... be
independent rv's with this distribution. Under these notations one obtains
from (9) by change of variable and monotonicity of $\Psi _{n}$
$$%
\Phi _{n+1}(s)=\Bbb E\Phi _{n}(s-\log \theta)-\Psi _{n}(s)\geq \Bbb E%
\Phi _{n}(s-\log \theta)-\Psi _{0}(s), n=1,2,...%
\tag 12
$$%
Consider the random walk $S_{0}=0$, $S_{n}=-\sum_{i=1}^{n}\log \theta _{i}$,
$n=1,2,...$ On iterating (12) one gets
$$%
\Phi _{n}(s)\geq \Bbb E\Phi _{0}(s+S_{n})-\Bbb E\sum_{i=0}^{n-1}\Psi
_{0}(s+S_{i}), n=1,2,...%
$$%
Note that $\Phi _{0}(s):=e^{s}(1-\exp (-e^{-s}))1\{s\geq \log s_{0}\}$ and $\Psi  _{0}(s):= $ \newline $e^{s}(e^{-s}-(1-\exp (-e^{-s})))1\{s\geq \log s_{0}\}$. Since $ \Bbb E\log \theta _{i}=$ \newline $=\lambda \int_{0}^{\infty }h(z)\log h(z)dz<0$ then
by the strong law of large numbers $S_{n}\rightarrow +\infty $ a.s., as $%
n\rightarrow \infty $. Consequently by dominated convergence
$$%
\underset{n\rightarrow \infty } \to {\lim }\Bbb E\Phi _{0}(s+S_{n})=m.%
\tag 13
$$%
Let us verify that $\Psi _{0}(s)$ is directly Riemann integrable (dRi). This
together with the renewal theorem on the whole line (Feller (1966), Theorem
XI.1, p. 368) will imply
$$%
\underset{s\rightarrow +\infty } \to {\lim }\Bbb E \sum_{i=0}^{\infty }\Psi
_{0}(s+S_{i})=0.%
\tag 14
$$%
Since $e^{-s}\Psi _{0}(s)$ is a decreasing function in $s$ and \newline $%
\int_{-\infty }^{\infty }\Psi
_{0}(s)ds=\int_{0}^{s_{0}}s^{-2}(s-(1-e^{-s}))ds<\infty $ (the latter
integrand is equivalent to a constant near $0$) then $\Psi _{0}(s)$ is
indeed dRi. Now according to (13), (14)
$$%
m\leq \underset{s\rightarrow +\infty } \to {\lim \inf } \underset{n\rightarrow
\infty } \to {\lim }\Phi _{n}(s)\leq \underset{s\rightarrow +0} \to {\lim \inf }%
s^{-1}(1-\varphi ^{(2)}(s))\overset{(11)} \to {\leq }\underset{s\rightarrow +0} \to {\lim \inf } s^{-1}(1-\varphi ^{(1)}(s))%
$$%
whence from (10) $\underset{s\rightarrow +0} \to {\lim }s^{-1}(1-\varphi
^{(1)}(s))=m$ or equivalently $\mu ^{(1)} =\mu ^{\ast }\in \Cal P_{h,m}^{+}$.

To prove uniqueness let us assume on the contrary that there exists another
LST $\widetilde{\varphi }(s)$ with $\underset{s\rightarrow +0} \to {\lim }%
s^{-1}(1-\widetilde{\varphi }(s))=m$ that satisfies (5). Choose in Lemma 3.1 \newline $\varphi _{1}:=\widetilde{\varphi }$, $\varphi _{2}:=\varphi ^{(1)}$, and for each positive integer $n$, $C_{n}:=\theta _{1}...\theta _{n}$. Note that
$C_{n}=\exp (-S_{n})\rightarrow 0$ a.s, as $n\rightarrow \infty $. By using
an elementary inequality \newline $\left| e^{-x}-e^{-y}\right| \leq \left| x-y\right|
,x,y\in \Bbb R$ one concludes from (5) $\psi (s)\leq \Bbb E%
\psi (\theta s)$, $s\geq 0$. On iterating the latter inequality $n$ times
one gets $\psi (s)\leq \Bbb E\psi (C_{n}s)$ which in turn implies that
Lemma 3.1 does apply. Hence  $\widetilde{\varphi }\equiv $ $\varphi ^{(1)}$
which completes the proof.%
\qed\enddemo%

Now we are ready to prove Proposition 2.1.

\demo {Proof of Proposition 2.1}%
We will prove both parts of the Proposition simultaneously. First of all note that the statement ''$\rho $ is a probability measure'' and condition (6) with $\lambda \int_{0}^{\infty}h(z)dz=1$ are equivalent.%
\newline
Further suppose that the conditions of part (b) are in force and the SNT $%
\Bbb T_{h,\lambda }$ has a fixed point $\mu ^{\ast }$. Hence $\varphi
^{\ast }(s)=\int_{0}^{\infty }e^{-sx}\mu ^{\ast }(dx)$ satisfies (5), that
is, $\varphi ^{\ast }(s)=$
$$%
=\exp\{-\lambda \int_{0}^{\infty }(1-\varphi ^{\ast }(sh(u)))du\}=\exp%
\{\lambda \int_{a}^{b}(1-\varphi ^{\ast }(sz))h^{\leftarrow }(dz)\}=%
$$%
$$%
=\exp\{-\lambda \int_{a}^{b}(1-\varphi ^{\ast}(sz)) z^{-1}\rho (dz)\},%
$$%
where $\rho (dz)=-\lambda zh^{\leftarrow}(dz)$ is a probability measure. The latter follows from (6) and (8).
In view of Lemma 3.2 condition (8) implies $m=\int_{0}^{\infty }x\mu ^{\ast
}(dx)<\infty $, for some $m>0$. W.l.o.g. we may and do assume $m=1$, and
therefore \newline $\underset{s\rightarrow +0} \to {\lim }s^{-1}(1-\varphi ^{\ast }(s))=1$.

On the other hand, assume that an rv $\eta $ with $\Bbb E\eta =1$
satisfies (2) and the conditions of part (a) hold. Then the LT $\varphi (s)=%
\Bbb E e^{-s\eta \text{ }}$solves
$$%
\varphi ^{^{\prime }}(s)=\varphi (s)\int_{0}^{\infty }\varphi ^{^{\prime
}}(sz) \rho (dz).%
$$%
Note that $-\varphi ^{^{\prime }}(s)$ is the LST of probability measure $\mu
_{sb}(dx):=x\mu (dx)$. By using Fubini's Theorem one has $\log \varphi
(s)=\int_{0}^{s}[\log \varphi (u)]^{^{\prime
}}du=\int_{0}^{s}\int_{0}^{\infty }\varphi ^{^{\prime }}(uz)\rho
(dz)du=\int_{0}^{\infty }z^{-1}\rho (dz)(\varphi (sz)-1)$ or equivalently
$$%
\varphi (s)=\exp \{-\int_{0}^{\infty }(1-\varphi (sz))z^{-1}\rho (dz)\}.%
$$%
Consider the (possibly $\sigma $-finite) measure $\nu (dz)=(\lambda
z)^{-1}\rho (dz)$ and the (right-continuous and non-increasing) function $%
V(x)=-\int_{x}^{b}\nu (dz)$. Now define the response function $h$ via its
generalized inverse $h^{\leftarrow }(x)=V(x)$. This implies $\rho
(dz)=-\lambda zh^{\leftarrow }(dz)$ or equivalently (7).

It remains to verify that $\varphi ^{\ast }(s)=\varphi (s)$ which can be
achieved in the same manner as in the proof of Lemma 3.2 (use Lemma 3.1).
The proof is completed.%
\qed\enddemo%

\demo {Proof of Proposition 1.1}%
Assume that $\Bbb P\{A=0\}=0$ and there exists non-zero $\mu $ solving (2), but
$\Bbb E\log A\geq 0$. By Theorem 1.6(a) of Vervaat (1979) the latter implies that $\Cal L(\eta
_{sb})=\delta _{c}$, for some $c\in \Bbb R$. This is possible if $c=0$,
the case excluded by us. A contradiction.\newline
Assume now that (3) holds and $\Bbb P\{A=0\}=0$. Consider the SNT $\Bbb T_{h,1}$ with the response function $h$ defined via its generalized inverse $h^{\leftarrow}(z)$ as in (7). Since $\Bbb E \log A<0$ implies $\int_{0}^{\infty}h(z)\log h(z)dz<0$, Lemma 3.2 applies that allows us to conclude that given $m\in (0,\infty )$ $\Bbb T_{h,1}$ has a unique fixed point $\mu ^{\ast }$ with $m=\int_{0}^{\infty }x\mu ^{\ast }(dx)$. It remains to appeal to Proposition 2.1(b). This finishes the proof of Proposition 1.1.%
\qed\enddemo%
To prove Proposition 1.2 we need the following result.
\proclaim {Lemma 3.3}%
Let $h$ be a positive measurable
function and for some $p>0$ \newline $\int_{0}^{\infty }h^{p}(u)du<\infty $. If $%
\Bbb E\xi _{1}^{p}<\infty $ then $\Bbb E(\sum_{i=1}^{\infty }\xi
_{i}h(\tau _{i}))^{p}<\infty $.%
\endproclaim%
\demo {Proof}%
Shot noise distributions are ID, and therefore as it is
well-known, \newline $\Bbb E(\sum_{i=1}^{\infty }\xi _{i}h(\tau _{i}))^{p}<\infty
$ iff $\int_{1}^{\infty }x^{p}M(dx)<\infty $ where $M$ is the L\'{e}vy
measure. Since $\int_{0}^{\infty }x^{p}M(dx)=\Bbb E \xi
_{1}^{p}\int_{0}^{\infty }h^{p}(u)du<\infty $, the assertion follows.%
\qed\enddemo%
\demo {Proof of Proposition 1.2(a)}%
1) Given the rv $A$ with $\Bbb P%
\{A=0\}=0$ and $\Bbb E A^{q-1}<1$, $q\in (1,2)$, define the SNT $\Bbb T%
_{h,1}$ with the response function $h$ whose generalized inverse $%
h^{\leftarrow }$ is defined by (7). By Proposition 2.1(a), first, $%
\int_{0}^{\infty }h(u)du=1$ from which one concludes that the restriction of $%
\Bbb T_{h,1}$ to the set $\Cal P_{m}^{+}(q)$ is well-defined and
for any $\theta \in \Cal P_{m}^{+}(q)$ $\int_{0}^{\infty }x(\Bbb T%
_{h,1}\theta )(dx)=m$; second, $\int_{0}^{\infty }h^{q}(u)du<1$ that in
conjunction with Lemma 3.3 imply $\int_{0}^{\infty }x^{q}(\Bbb T%
_{h,1}\theta )(dx)<\infty $. Thus $\Bbb T_{h,1}$ maps $\Cal P%
_{m}^{+}(q)$ into itself. The key observation for what follows is that for
any $\theta \in \Cal P_{m}^{+}(q)$ and $\rho =\Cal L (A)$ where $A$
satisfies the above conditions, the map $\Bbb U_{\rho \text{ }}$ is
well-defined and moreover
$$%
\Bbb U_{\rho \text{ }}\theta =\Bbb T_{h,1}\theta.%
$$%
Formally this can be verified by recording the corresponding LST's. To prove
1), it remains to add that according to Lemma 3.5 of Iksanov, Jurek (2002)
the restriction of $\Bbb T_{h,1}$ defined on a complete metric space $(%
\Cal P_{m}^{+}(q),r_{q})$ is a strict contraction, and apply the Banach
Fixed Point Theorem. \newline
2) The preceding display and Proposition 2.1(b) imply that $\mu $ solves (2). Assume now (4) holds and at once note that the assertion $\Bbb E \eta
=m$ is obvious. We intend to verify that (4) implies $\Bbb E\eta
^{1+p}<\infty $. W.l.o.g. set $m=1$ and consider two cases:\newline
1) $p\in \Bbb N$: if $p=1$ then (2) implies $\Bbb E\eta _{sb}=%
\Bbb E\eta /(1-\Bbb E A)<\infty $ since by the definition of the
size-biased distribution $\Bbb E \eta <\infty $. Further we proceed by induction. Assume it is already known that for some $%
1<k<p$, $k\in \Bbb N$
$$%
\Bbb E \eta ^{1+k}=\Bbb E\eta _{sb}^{k}<\infty .%
$$%
By assumption (4) the function $g(x)$ is log-convex on $(0,p)$ which implies
$g(l)=\Bbb E A^{l}<1$ for all $l\in (0,p)$. Therefore $\Bbb E%
A^{k+1}<1 $. Now an appeal to (2) allows us to write formally
$$%
\gather%
\Bbb E\eta ^{2+k}=\Bbb E\eta _{sb}^{1+k}=%
\Bbb E (A\eta _{sb}+\eta )^{1+k}=%
\\%
=\sum_{i=0}^{1+k} \binom {1+k}i%
\Bbb E A^{i}\Bbb E\eta ^{i+1}\Bbb E \eta
^{1+k-i}=%
\\%
=\Bbb E A^{1+k}\Bbb E \eta _{sb}^{1+k}+\sum_{i=0}^{k}\binom {1+k}i%
\Bbb E A^{i}\Bbb E \eta ^{i+1}\Bbb E\eta%
^{1+k-i}:=%
\\%
:=\Bbb E A^{1+k}\Bbb E\eta ^{2+k}+t_{k},%
\endgather%
$$%
which yields $\Bbb{E}\eta ^{2+k}=t_{k}/(1-\Bbb E A^{1+k})<\infty $ by
noting that $t_{k}$ consists of the finite number of finite terms. This
completes the study of this case.\newline
2) $p\notin \Bbb N $: denote by $\alpha $ the fractional part of $p$ then
there exists $n\in \Bbb N\cup \{0\}$ such that $p=n+\alpha $. If $n=0$
the assertion follows from the first part of the Proposition where in fact
was shown that if $\Bbb E A^{\alpha }<1$ then $\mu =\Cal L (\eta )\in
\Cal P_{m}^{+}(1+\alpha )$. Suppose we have already verified that for
some $0<k<n$, $k\in \Bbb N$
$$%
\text{\ }\Bbb E\eta ^{1+k+\alpha }=\Bbb E\eta _{sb}^{k+\alpha
}<\infty.%
\tag 15
$$%
Set $k_{\alpha }:=1+k+\alpha $. Now (2) implies that one may write formally
$$%
\gather%
(\Bbb E (\eta ^{k_{\alpha}+1}))^{k_{\alpha}^{-1}}=(\Bbb E \eta
_{sb}^{k_{\alpha}})^{k_{\alpha}^{-1}}=%
\\%
=(\Bbb E (A\eta _{sb}+\eta)^{k_{\alpha}})^{k_{\alpha}^{-1}}\leq%
\\%
\leq (\Bbb E A^{k_{\alpha}})^{k_{\alpha}^{-1}}(\Bbb E \eta ^{k_{\alpha
}+1})^{k_{\alpha}^{-1}}+(\Bbb E \eta ^{k_{\alpha}})^{k_{\alpha}^{-1}},%
\endgather%
$$%
the inequality being implied by the triangle inequality in the space $%
L_{k_{\alpha }}$. Now the latter inequality can be rewritten as follows
$$%
\Bbb E \eta ^{k_{\alpha }+1}\leq \Bbb E \eta ^{k_{\alpha }}/(1-(%
\Bbb E A^{k_{\alpha }})^{k_{\alpha }^{-1}})^{k_{\alpha }}<\infty \,%
$$%
where finiteness is implied by (15) and log-convexity of the function $g(x)$
that guarantees $1-(\Bbb E A^{k_{\alpha }})^{k_{\alpha }^{-1}}>0$. A
usual inductive argument completes the proof.%
\qed\enddemo%
\demo {Proof of Proposition 1.2(b)}%
While for $p>1$ we use an elementary
inequality \newline $(x+y)^{p}\geq x^{p}+y^{p}$, for $x$, $y\geq 0,$ to obtain $
\Bbb E \eta _{sb}^{p}=\Bbb E (\eta +A\eta _{sb})^{p}\geq \Bbb E%
\eta ^{p}+\Bbb E A^{p}\Bbb E \eta _{sb}^{p}>\Bbb E A^{p}\Bbb E %
\eta _{sb}^{p}$; if $p\in (0,1)$, then by a variant of Minkowski inequality $%
(\Bbb E \eta _{sb}^{p})^{1/p}=(\Bbb E (\eta +A\eta
_{sb})^{p})^{1/p}\geq (\Bbb E \eta ^{p})^{1/p}+(\Bbb E A^{p})^{1/p}(%
\Bbb E \eta _{sb}^{p})^{1/p}>(\Bbb E A^{p})^{1/p}(\Bbb E \eta
_{sb}^{p})^{1/p}$. Either of them implies the desired. The case $p=1$
is trivial.%
\qed\enddemo%
\demo {Proof of Proposition 1.3}%
Clearly, the proof could be done by using
the connection to the shot noise distributions (Proposition 2.1). However, we
choose to use a more direct way. Let $\omega $ be a non-negative ID
distribution with the L\'{e}vy measure $N$ and shift $\gamma \geq 0$. The
following representation due to Steutel (1970, p.86) is well-known
$$%
\int_{0}^{x}y\omega (dy)=\int_{0+}^{x}\omega [0,x-y]yN(dy)+\gamma
\omega [0,x].%
\tag 16
$$%
Note that $\gamma $ is merely a mass of an atom at zero of the measure $%
R(dx):=xN(dx).$ On the other hand, each distribution $\omega $ satisfying
(16) with appropriate $N$ and $\gamma $ is ID.

Putting $m:=1$ let us now rewrite equality (2) in terms of distribution
functions to obtain
$$%
\mu _{sb}[0,x]=\int_{0}^{x}y\mu (dy)=\int_{0}^{x}\mu [0,x-y]yM(dy),%
\tag 17
$$%
where the measure $M$ satisfies the relation
$$%
\int_{0}^{x}yM(dy)=\int_{a}^{b}\mu _{sb}[0,x/y]\rho (dy),%
\tag 18
$$%
or equivalently $M(dx)=x^{-1}\Cal L (A\eta _{sb})(dx)$. By comparing
(18), (17) with (16) one immediately concludes that $\mu $ is ID with the L\'{e}vy measure $M$. Its shift is $0$ because $\Bbb P \{A=0\}=0$ which
implies that the measure $S(dx):=xM(dx)$ is atomless at zero. To complete
the proof recall that finiteness of the L\'{e}vy measure is necessary and
sufficient condition for a distribution to be compound Poisson. In our case
this requirement reduces to the integrability of the measure $x^{-1}\rho
(dx) $ at the origin.%
\qed\enddemo%
\demo {Proof of Proposition 1.4}%
The part (a) follows from Theorem
1.1(b) of Iksanov, Jurek (2002) taking into account Proposition 2.1. Turn
now to the proof of part (b). W.l.o.g. we only consider the
case $m=1$. By Proposition 1.3 $\mu $ is ID with the L\'{e}vy measure $%
M(dx)=x^{-1}\Cal L (A\eta _{sb})(dx)$. By Corollary 25.8 of Sato (1999)
$$%
\Bbb E e^{s\eta _{sb}}=\int_{0}^{\infty }e^{sx}x\mu (dx)<\infty \text{ \
iff }\int_{1}^{\infty }e^{sx}xM(dx)<\infty.%
\tag 19
$$%
Put $c=$esssup $\rho \in (0,1)$. It suffices to verify that if $\Bbb E%
e^{s\eta _{sb}}<\infty $ then $\Bbb E e^{(s/c)\eta _{sb}}<\infty $ as
well. In view of part (a), there exists $s_{0}=s_{0}(c)>0$ such that $\Bbb %
E\exp (s_{0}c\eta _{sb})<\infty $. Then $\Bbb E \exp (s_{0}A\eta
_{sb})<\infty $ and consequently $\int_{1}^{\infty }e^{s_{0}x}xM(dx)<\infty $%
. It remains to use (19) which finishes the proof.%
\qed\enddemo%
\demo {Proof of Proposition 1.5}%
Suppose $A\in (0,1]$ a.s. Then $%
\Bbb E A^{p}<1$, for all $p>0$. Consequently, by Proposition 1.2 (a) $%
\Bbb E \eta ^{p}<\infty $, for all $p>0$. In the sequel ''$\rightarrow $'' means ''uniquely determines''. Let us imagine the chain and prove its validity $\{%
\Bbb E \eta ^{n}\}_{n\in \Bbb N} \rightarrow \{\Bbb E %
A^{n}\}_{n\in \Cal N}\rightarrow \Cal L (A)\rightarrow \Cal L%
(\eta )$, where $\eta $ is the rv with fixed mean $m<\infty $. \newline 1) $\{%
\Bbb E \eta ^{n}\}_{n\in \Bbb N}\rightarrow \{\Bbb E%
A^{n}\}_{n\in \Bbb N}$ via the relation
$$%
\Bbb E \eta ^{n+1}=\sum_{k=0}^{n} \binom n k \Bbb E A^{k}\Bbb E \eta ^{k+1}\Bbb E \eta ^{n-k}.%
$$%
2) $\{\Bbb E A^{n}\}_{n\in \Bbb N}\rightarrow \Cal L (A)$ since $%
A\in (0,1]$ a.s. (the moments problem). \newline 3) $\Cal L (A)\rightarrow
\Cal L (\eta )$ which is obvious by Proposition 1.1(a)(or 1.2 (a)). \newline
Conversely, if $\Cal L (\eta )$ is completely determined by its moments
then $\Bbb E \eta ^{p}<\infty $, for all $p>0$ which in view of Theorem
1.2(b) implies $\Bbb E A^{p}<1$, for all $p>0$. Clearly, this makes $A\in
(0,1]$ a.s. and completes the proof.%
\qed\enddemo%

\enddocument